\title{Characteristically simple Beauville groups, II:\\
low rank and sporadic groups}

\author{Gareth A. Jones\\
School of Mathematics\\
University of Southampton\\
Southampton SO17  1BJ, U.K.\\
{\tt G.A.Jones@maths.soton.ac.uk}
}

\documentclass[12pt]{article}
\usepackage{color}
\usepackage[latin1]{inputenc}
\usepackage{a4wide}
\usepackage{latexsym}
\usepackage{amsmath}
\usepackage{amssymb}
\usepackage{amsfonts}
\newtheorem{thm}{Theorem}[section]
\newtheorem{lemma}[thm]{Lemma}
\newtheorem{cor}[thm]{Corollary}
\newtheorem{prop}[thm]{Proposition}
\date{}
\begin{document}

\maketitle

\begin{abstract}
A Beauville surface is a rigid complex surface of general type, isogenous to a higher product by the free action of a finite group, called a Beauville group. Here we consider which characteristically simple groups can be Beauville groups. We show that if $G$ is a cartesian power of a simple group $L_2(q)$, $L_3(q)$, $U_3(q)$, $Sz(2^e)$, $R(3^e)$, or of a sporadic simple group, then $G$ is a Beauville group if and only if it has two generators and is not isomorphic to $A_5$.

\end{abstract}

 (secondary).

\section{Introduction}

A Beauville surface $\cal S$ of unmixed type is a complex algebraic surface which is isogenous to a higher product (that is, ${\cal S}=({\cal C}_1\times{\cal C}_2)/G$, where ${\cal C}_1$ and ${\cal C}_2$ are complex algebraic curves of genus $g_i\ge 2$, and $G$ is a finite group acting freely on their product), and is rigid in the sense that $G$ preserves the factors ${\cal C}_i$, with ${\cal C}_i/G\cong {\mathbb P}^1({\mathbb C})$ ($i=1, 2$) and the induced covering $\beta_i:{\cal C}_i\to{\mathbb P}^1({\mathbb C})$ ramified over three points.
(In this paper we will not consider Beauville surfaces of mixed type, where elements of $G$ transpose the factors ${\cal C}_i$.) The first examples, with ${\cal C}_1={\cal C}_2$ Fermat curves, were introduced by Beauville~\cite[p.~159]{Bea} in 1978. Subsequently, these surfaces have been intensively studied by geometers such as Bauer, Catanese and Grunewald (see~\cite{BCG05, BCG06, Cat}, for instance).

Recently, there has been conderable interest in determining which groups $G$, known as Beauville groups, can arise in this way. It is easily shown that the alternating group $A_5$ is not a Beauville group.  In 2005 Bauer, Catanese and Grunewald~\cite{BCG05} conjectured that every other non-abelian finite simple group is a Beauville group. In 2010, more or less simultaneously, Garion, Larsen and Lubotzky~\cite{GLL} proved that this is true with finitely many possible exceptions, while Guralnick and Malle~\cite{GM} and Fairbairn, Magaard and Parker~\cite{FMP} proved that it is completely true. It is natural to consider other classes of finite groups, and this paper is part of a project to extend this result to characteristically simple groups, those with no characteristic subgroups other than $1$ and $G$.

A finite group is characteristically simple if and only if it is isomorphic to a cartesian power $H^k$ of a simple group $H$. In~\cite{Cat}, Catanese showed that the abelian Beauville groups are those isomorphic to $C_n^2$ for some $n$ coprime to $6$; these are characteristically simple if and only if $n$ is prime, so we may assume that $H$ is a non-abelian finite simple group. 

The ramification condition on the coverings $\beta_i$ implies that a Beauville group is a quotient of a triangle group (in two ways), so it is a $2$-generator group. Clearly we require that $G\not\cong A_5$, and in~\cite{Jon13} it was conjectured that these two conditions are also sufficient:

\medskip

\noindent{\bf Conjecture} {\it Let $G$ be a non-abelian finite characteristically simple group. Then $G$ is a Beauville group if and only if it is a $2$-generator group not isomorphic to $A_5$.}

\medskip

This was proved in~\cite{Jon13} in the case where $G=H^k$ with $H$ an alternating group $A_n$. (Further details on the background to this problem are given there, while~\cite{Jon12} gives a more general survey of Beauville groups.) By the classification of finite simple groups, this leaves the simple groups of Lie type, together with the $26$ sporadic simple groups. Our aim here is to verify the conjecture for five families of groups of small Lie rank, specifically the projective special linear groups $L_2(q)$ and $L_3(q)$, the unitary groups $U_3(q)$, the Suzuki groups $Sz(2^e)$ and the `small' Ree groups $R(3^e)$, together with the sporadic simple groups (see~\cite{ATLAS, Wil} for notation and properties of these and other finite simple groups). 

The $2$-generator condition on $H^k$ can be restated as follows. For any $2$-generator finite group $H$, there is an integer $c_2(H)$ such that $H^k$ is a $2$-generator group if and only if $k\le c_2(H)$. If $H$ is a non-abelian finite simple group (and hence a $2$-generator group, by the classification), then $c_2(H)$ is equal to the number $d_2(H)$ of orbits of ${\rm Aut}\,H$ on ordered pairs of generators of $H$. Thus in this case, $H^k$ is a $2$-generator group if and only if $k\le d_2(H)$. 

For the groups $L_2(q)=SL_2(q)/\{\pm I\}$, which are simple for all prime powers $q\ge 4$, with $L_2(4)\cong L_2(5)\cong A_5$, we will prove:

\begin{thm}\label{mainthm1}
Let $G=H^k$, where $k\ge 1$ and $H=L_2(q)$ for a prime power $q\ge 4$. Then the following are equivalent:
\begin{description}
\item {\rm(a)} $G$ is a Beauville group;
\item {\rm(b)} $G$ is a $2$-generator group not isomorphic to $A_5$;
\item {\rm(c)} $k\le d_2(H)$, and if $q\le 5$ then $k>1$.
\end{description}
\end{thm}

The groups $L_3(q)$ and $U_3(q)$ are simple for all prime powers $q\ge 2$ and $3$ respectively, while the Suzuki groups $Sz(2^e)$ and the Ree groups $R(3^e)$ are simple for all odd $e\ge 3$. For these and for the sporadic simple groups, we have a slightly simpler formulation:

\begin{thm}\label{mainthm2}
Let $G=H^k$, where $k\ge 1$ and $H$ is a simple group $L_3(q)$, $U_3(q)$, $Sz(2^e)$ or $R(3^e)$,  or a sporadic simple group. Then the following are equivalent:
\begin{description}
\item {\rm(a)} $G$ is a Beauville group;
\item {\rm(b)} $G$ is a $2$-generator group;
\item {\rm(c)} $k\le d_2(H)$.
\end{description}
\end{thm}

In Theorems~\ref{mainthm1} and~\ref{mainthm2}, the implications (a) $\Rightarrow$ (b) $\Rightarrow$ (c) are straightforward (see Sections~\ref{Bmp} and~\ref{Gcp}), and the main task is to prove that (c) implies (a). The case where $k=1$ having been dealt with by others (see above), it is thus sufficient to prove the following:

\begin{thm}\label{mainthm}
If $H$ is a simple group $L_2(q)$, $L_3(q)$, $U_3(q)$, $Sz(2^e)$ or $R(3^e)$,  or a sporadic simple group, then $H^k$ is a Beauville group for each $k=2,\ldots, d_2(H)$.
\end{thm}

It is hoped that the methods developed here and in~\cite{Jon13} may allow this result to be extended to all finite simple groups of Lie type, thus proving the above conjecture.
 
The results proved here and in~\cite{Jon13} yield further examples of Beauville groups, as follows. Let us say that two groups are {\sl orthogonal\/} if they have no non-identity epimorphic images in common. It is straightforward to show that a cartesian product of finitely many mutually orthogonal Beauville groups is a Beauville group: orthogonality ensures that the obvious triples generate the group. It follows that if $G_1,\ldots, G_r$ are characteristically simple Beauville groups, which are cartesian powers of mutually non-isomorphic simple groups, then $G_1\times\cdots\times G_r$ is also a Beauville group.

\section{Background and method of proof}\label{Bmp}

The proof of Theorem~\ref{mainthm} is based on the characterisation by Bauer, Catanese and Grunewald \cite{BCG05, BCG06} of Beauville groups as those finite groups $G$ such that
\begin{itemize}
\item $G$ is a smooth quotient of two hyperbolic triangle groups
\[\Delta_i=\Delta_i(l_i, m_i, n_i)=\langle A_i, B_i, C_i\mid A_i^{l_i}=B_i^{m_i}=C_i^{n_i}=A_iB_iC_i=1\rangle\]
for $i=1, 2$,
\item $\Sigma_1\cap\Sigma_2=1$, where $\Sigma_i$ is the set of conjugates of powers of the elements $a_i, b_i, c_i$ of $G$ which are images of the canonical generators $A_i, B_i, C_i$ of $\Delta_i$.
\end{itemize}
Here `smooth' means that $a_i, b_i$ and $c_i$ have orders $l_i, m_i$ and $n_i$, so that the kernel $K_i$ of the epimorphism $\Delta_i\to G$ is torsion-free, and thus a surface group. `Hyperbolic' means that $l_i^{-1}+m_i^{-1}+n_i^{-1}<1$, so that $\Delta_i$ acts on the hyperbolic plane $\mathbb H$. The isomorphism $G\cong\Delta_i/K_i$ induces an action of $G$ as a group of automorphisms of the compact Riemann surface (or complex algebraic curve) ${\cal C}_i={\mathbb H}/K_i$ of genus $g_i>1$, with the natural projection ${\cal C}_i\to{\cal C}_i/G\cong{\mathbb H}/\Delta_i\cong{\mathbb P}^1({\mathbb C})$ branched over three points (corresponding to the fixed points of $A_i, B_i$ and $C_i$ in $\mathbb H$). Since $\Sigma_i$ is the set of elements of $G$ with fixed points in ${\cal C}_i$, the condition $\Sigma_1\cap\Sigma_2=1$ is necessary and sufficient for $G$ to act freely on ${\cal C}_1\times{\cal C}_2$, so that ${\cal S}:=({\cal C}_1\times{\cal C}_2)/G$ is a complex surface. The generating triples $(a_i, b_i, c_i)$ for $G$ (always satisfying $a_ib_ic_i=1$) are said to have {\sl type\/} $(l_i, m_i, n_i)$, and the pair of them form a {\sl Beauville structure\/} of type $(l_1, m_1, n_1; l_2, m_2, n_2)$ in $G$.

Since $\Sigma_i$ is closed under taking powers, the condition $\Sigma_1\cap\Sigma_2=1$ is equivalent to the condition that
\[\Sigma_1^{(p)}\cap\Sigma_2^{(p)}=\emptyset\]
for every prime $p$ dividing $l_1m_1n_1$ and $l_2m_2n_2$, where $\Sigma_i^{(p)}$ denotes the set of elements of order $p$ in $\Sigma_i$. In many cases, this condition is simpler to verify.

If $G=H^k$ for some $H$ then the members of any generating triple for $G$ are $k$-tuples $a_i=(a_{ij})$, $b_i=(b_{ij})$, $c_i=(c_{ij})$ such that $(a_{ij}, b_{ij}, c_{ij})$ is a generating triple for $H$ for each $j=1,\ldots, k$, and these $k$ triples are mutually inequivalent under the action of ${\rm Aut}\,H$. If $H$ is a non-abelian finite simple group then conversely any $k$ mutually inequivalent generating triples for $H$ yield a generating triple for $G$. Moreover, if $G\not\cong A_5$ this triple is hyperbolic. To prove that $G$ is a Beauville group it is therefore sufficient to find two $k$-tuples of mutually inequivalent generating triples for $H$, and then to show that the resulting pair of generating triples for $G$ satisfy $\Sigma_1^{(p)}\cap\Sigma_2^{(p)}=\emptyset$ for all primes $p$ dividing $|H|$.

There are $k$-tuples of mutually inequivalent generating triples for $H$ if and only if $k\le d_2(H)$. In~\cite{Jon13}, suitable generating triples for the simple groups $H=A_n$ were exhibited by explicitly defining their members $a_{ij}, b_{ij}$ and $c_{ij}$, and then applying results on groups containing such permutations to show that they generate $H$. Here we adopt a similar approach for the groups $H=L_2(q)$, using Dickson's description of their maximal subgroups~\cite[Ch.~XII]{Dic} to show that various specific triples generate $H$. However, for the remaining simple groups it is more convenient to use a character-theoretic formula of Frobenius (see Section~\ref{Ct}) to prove the existence of triples of various types, and then again to use knowledge of the maximal subgroups of $H$ to show that they generate $H$.

In each case, the condition that $\Sigma_1^{(p)}\cap\Sigma_2^{(p)}=\emptyset$ for all $p$ is ensured by carefully choosing the two $k$-tuples of generating triples for $H$ so that powers of $a_1, b_1$ or $c_1$ of order $p$ have different supports in $\{1, 2, \ldots, k\}$ from those of $a_2, b_2$ or $c_2$. Several results which guarantee this were proved in~\cite{Jon13}, and for completeness their proofs are outlined in Section~\ref{Bscp}. These results are then applied to the groups $L_2(q)$ in Section~\ref{L2q}, and to the other families of simple groups in Sections~\ref{Suz} to~\ref{L3qU3q}. However, before doing this we will briefly discuss in Section~\ref{Gcp} the values of $k$ for which $H^k$ is a $2$-generator group.

\section{Generating cartesian powers}\label{Gcp}

If $H$ is a finite group, let
\[c_2(H)=\max\{k\in{\mathbb N}\mid H^k\; \hbox{is a $2$-generator group}\}.\]
If $H^k$ is a Beauville group then it is a $2$-generator group, and thus a quotient of the free group $F_2$ of rank $2$, so
\[k \le c_2(H)\le d_2(H)\]
where $d_2(H)$ denotes the number of normal subgroups $N$ of $F_2$ with $F_2/N\cong H$. Following Hall~\cite{Hal}, let $\phi_2(H)$ be the number of $2$-bases (ordered generating pairs) for $H$. Any $2$-base for $H$ determines an epimorphism $\theta: F_2\to H$, and hence a normal subgroup $N=\ker\theta$ of $F_2$ with $F_2/N\cong H$; conversely, every such normal subgroup arises in this way, with two $2$-bases corresponding to the same normal subgroup $N$ if and only if they are equivalent under an automorphism of $H$. Thus $d_2(H)$ is equal to the number of orbits of ${\rm Aut}\,H$ on $2$-bases for $H$. Since only the identity automorphism can fix a $2$-base, this action is semiregular, so we have the following:

\begin{lemma}
If $H$ is a finite group then
\[d_2(H)=\frac{\phi_2(H)}{|{\rm Aut}\,H|}.\]
\end{lemma}

If $H$ is a non-abelian finite simple group, and $F_2$ has $k$ normal subgroups with quotient $H$, then their intersection is a normal subgroup with quotient $H^k$. It follows then that $c_2(H)=d_2(H)$ for such a group $H$, so we have the following:

\begin{cor}
If $H$ is a non-abelian finite simple group and $H^k$ is a Beauville group, then
\[k\le c_2(H)=d_2(H)=\frac{\phi_2(H)}{|{\rm Aut}\,H|}.\]
\end{cor}

Since generating triples $(a, b, c)$ for a group correspond bijectively to its $2$-bases $(a, b)$, we also obtain the following useful characterisation of generating triples for $H^k$:

\begin{cor}\label{carttriples}
Let $H$ be a non-abelian finite simple group. Then $k$-tuples $a=(a_j)$, $b=(b_j)$ and $c=(c_j)$ form a generating triple for $H^k$ if and only if their components $(a_j, b_j, c_j)$ for $j=1,\ldots, k$ form $k$ mutually inequivalent generating triples for $H$. \hfill$\square$
\end{cor}

By results of Dixon~\cite{Dix}, Kantor and Lubotzky~\cite{KL}, and Liebeck and Shalev~\cite{LS},  a randomly-chosen pair of elements of a non-abelian finite simple group $H$ generate the whole group with probability approaching $1$ as $|H|\to\infty$, so
 \[d_2(H)\sim\frac{|H|^2}{|{\rm Inn}\,H|.|{\rm Out}\,H|}=\frac{|H|}{|{\rm Out}\,H|}\quad\hbox{as}\quad |H|\to\infty.\]
For each infinite family of such groups, $|{\rm Out}\,H|$ is constant or grows much more slowly than $|H|$, so $d_2(H)$ grows almost as fast as $|H|$. In particular, as $q=p^e\to\infty$ with $p$ prime, the infinite families considered here satisfy

\begin{itemize}
\item $d_2(L_2(q))\sim q^3/d^2e$ while $|L_2(q)|\sim q^3/d$, where $d=\gcd(2, q-1)$,
\item $d_2(L_3(q))\sim q^8/2d^2e$ while $|L_3(q)|\sim q^8/d$, where $d=\gcd(3, q-1)$,
\item $d_2(U_3(q))\sim q^7/2d^2e$ while $|U_3(q)|\sim q^7/d$, where $d=\gcd(3, q+1)$,
\item $d_2(Sz(q))\sim q^5/e$ while $|Sz(q)|\sim q^5$,
\item $d_2(R(q))\sim q^7/e$ while $|R(q)|\sim q^7$.
\end{itemize}

One can illustrate how close the asymptotic estimate $|H|/|{\rm Out}\,H|$ is to $d_2(H)$ as follows. For any finite group $H$ we have
\[\phi_2(H)=|H^2\setminus \bigcup_MM^2|\ge |H|^2-\sum_M|M|^2,\]
where $M$ ranges over the maximal subgroups of $H$. If $H$ is perfect then each such $M$ has $|H:M|$ conjugates, so
\[\phi_2(H)\ge |H|^2\left(1-\sum_{i=1}^r\frac{1}{|H:M_i|}\right),\]
where $M_i$ ranges over a set of representatives of the $r$ conjugacy classes of maximal subgroups of $H$. It follows that if $H$ also has trivial centre (and in particular if $H$ is a non-abelian finite simple group) then
\[1\ge \frac{d_2(H)}{|H|/|{\rm Out}\,H|}\ge 1-\sum_{i=1}^r\frac{1}{|H:M_i|}.\]
When $H$ is large, the sum of the right is typically very small. Thus, for the simple groups $H=L_2(q)$, $U_3(q)$, $Sz(q)$ and $R(q)$ we have
\[\sum_{i=1}^r\frac{1}{|H:M_i|}\sim\frac{1}{q^m}\quad\hbox{as}\quad q\to\infty,
\]
where $m=1, 3, 2$ and $3$ respectively, with the sum dominated by the term corresponding to the doubly transitive permutation representation of $H$ of degree $q^m+1$; for $H=L_n(q)$ with $n\ge 3$ the sum grows like $2/q^{n-1}$, corresponding to the two doubly transitive representations of degree $(q^n-1)/(q-1)$, on the points and hyperplanes of ${\mathbb P}^{n-1}(q)$.

The situation is similar for the sporadic simple groups: for instance, if $H$ is O'Nan's group $O$'$N$, with $r=13$ conjugacy classes of maximal subgroups, we have
\[\frac{|H|}{|{\rm Out}\,H|}=\frac{460815505920}{2}=230407752960\]
and
\[\sum_{i=1}^r\frac{1}{|H:M_i|}=0.00001726863378\ldots,\]
so
\[230407752960\ge d_2(H)\ge 230403774132.\]

\section{Beauville structures in cartesian powers}\label{Bscp}

Corollary~\ref{carttriples} describes the generating triples in a cartesian power $G=H^k$ of a non-abelian finite simple group $H$, in terms of those for $H$. To avoid confusion, we will use notations such as ${\mathcal T}=(a, b, c)$ and $T=(x, y, z)$ for generating triples in $G$ and $H$ respectively. We now consider sufficient conditions for two triples ${\mathcal T}_i=(a_i, b_i, c_i)$ in $G$ to satisfy $\Sigma_1^{(p)}\cap\Sigma_2^{(p)}=\emptyset$ for all primes $p$, so that they form a Beauville structure for $G$.

If $g=(g_j)\in\Sigma_i^{(p)}$ for some prime $p$, then $g$ has order $p$ and is conjugate in $G$ to a power of some $d_i=a_i$, $b_i$ or $c_i$, so for each $j\in{\mathbb N}_k:=\{1,\ldots, k\}$ its $j$th coordinate $g_j$ has order $1$ or $p$. The support
\[{\rm supp}(g)=\{j\in{\mathbb N}_k\mid g_j\ne 1\}\]
of $g$ is then the $p$-{\sl summit\/} $S_p(d_i)$ of $d_i$, defined to be the set of $j\in{\mathbb N}_k$ for which the power of $p$ dividing the order of its $j$th coordinate $d_{ij}$ is greatest. The $p$-summit $S_p({\mathcal T})$ of a triple ${\mathcal T}=(a, b, c)$ in $G$ is the set $\{S_p(a), S_p(b), S_p(c)\}$ of subsets of ${\mathbb N}_k$. The following is obvious:

\begin{lemma}\label{disjointsummits}
Let $H$ be any finite group. If two triples ${\mathcal T}_i=(a_i, b_i, c_i)$ ($i=1, 2$) in $G=H^k$ have disjoint $p$-summits $S_p({\mathcal T}_i)$ for some prime $p$, then $\Sigma_1^{(p)}\cap\Sigma_2^{(p)}=\emptyset$. \hfill$\square$
\end{lemma}

An element of $H$ is $p$-{\sl full\/} if its order is divisible by the highest power of $p$ dividing the exponent ${\rm exp}(H)$ of $H$. If $d_i$ has at least one $p$-full coordinate $d_{ij}$, then
\[S_p(d_i)=F_p(d_i):=\{j\in{\mathbb N}_k\mid d_{ij}\; \hbox{is $p$-full}\}.\]
If each element $d_i$ of a triple ${\mathcal T}_i$ in $G$ has a $p$-full coordinate, this makes it easier to determine $S_p({\mathcal T}_i)$, and hence to ensure that pairs of such triples ${\mathcal T}_i$ satisfy Lemma~\ref{disjointsummits}.

Given a triple $T=(x, y, z)$ in $H$, let $\nu_p(T)$ be the number of $p$-full elements among $x, y$ and $z$. Two triples $T_i$ ($i=1, 2$) in $H$ are $p$-{\sl distinguishing} if $\nu_p(T_1)\ne\nu_p(T_2)$, and {\sl strongly} $p$-{\sl distinguishing} if, in addition, whenever $\nu_p(T_i)=0$ then either $p^2$ does not divide $\exp(H)$ or $p$ does not divide any of the three periods of $T_i$.

\begin{lemma}\label{stronglypdist}
Suppose that a non-abelian finite simple group $H$ has a set $\{(T_{1,s}, T_{2,s})\mid s=1,\ldots, t\}$ of ordered pairs $(T_{1,s}, T_{2,s})$ of generating triples such that
\begin{enumerate}
\item for each prime $p$ dividing $|H|$ there is some $s=s(p)\in\{1,\ldots, t\}$ such that $T_{1,s}$ and $T_{2,s}$ are a strongly $p$-distinguishing pair;
\item for each $i=1, 2$ the $3t$ triples consisting of $T_{i,1},\ldots, T_{i,t}$ and their cyclic permutations are mutually inequivalent.
\end{enumerate}
Then $G:=H^k$ is a Beauville group for each $k=3t,\ldots, d_2(H)$.
\end{lemma}

\noindent{\sl Proof.} Let $T_{i,s}=(x_{i,s}, y_{i,s}, z_{i,s})$ for each $i=1, 2$ and $s=1,\ldots., t$. Since $3t\le k\le d_2(H)$, for each $i$ one can form a generating triple ${\mathcal T}_i=(a_i, b_i, c_i)$ for $G$ by using $k$ inequivalent generating triples for $H$ in the different coordinate positions, with $T_{i,s}$ and its two cyclic permutations in the $j$th positions where $j=3s-2$, $3s-1$ or $3s$ for $s=1,\ldots, t$.

Suppose that $g\in\Sigma_1^{(p)}\cap\Sigma_2^{(p)}$ for some prime $p$. There is a strongly $p$-distinguishing pair $(T_{1,s}, T_{2,s})$, with $\nu_p(T_{1,s})\ne\nu_p(T_{2,s})$. First suppose that $\nu_p(T_{i,s})>0$ for $i=1, 2$, so that each element $d_i=a_i, b_i$ or $c_i$ of ${\mathcal T}_i$ has at least one $p$-full coordinate, and hence $S_p(d_i)=F_p(d_i)$. It follows that
\begin{equation}
|{\rm supp}(g)\cap\{3s-2, 3s-1, 3s\}|=\nu_p(T_{i,s})
\end{equation}
for $i=1$ and $2$, which is impossible since $\nu_p(T_{1,s})\ne\nu_p(T_{2,s})$. Otherwise, we may suppose without loss of generality that $\nu_p(T_{1,s})>0=\nu_p(T_{2,s})$, with $(1)$ satisfied for $i=1$ but not for $i=2$. Some $d_2=a_2, b_2$ or $c_2$ must then have a non-$p$-full coordinate of order divisible by $p$ in position $3s-2, 3s-1$ or $3s$, which is impossible if the periods of $T_{2,s}$ are coprime to $p$ or if $p^2$ does not divide ${\rm exp}(H)$, as we are assuming. Thus $\Sigma_1^{(p)}\cap\Sigma_2^{(p)}=\emptyset$, so the triples ${\mathcal T}_i$ form a Beauville structure for $G$. \hfill$\square$



\medskip

If we simply assume that the pairs $T_{1,s}$ and $T_{2,s}$ are $p$-distinguishing, we have the following rather weaker conclusion:

\begin{lemma}\label{pdist}
Suppose that a non-abelian finite simple group $H$ has a set ${\cal T}=\{(T_{1,s}, T_{2,s})\mid s=1,\ldots, t\}$ of ordered pairs $(T_{1,s}, T_{2,s})$ of generating triples such that
\begin{enumerate}
\item for each prime $p$ dividing $|H|$ there is some $s=s(p)\in\{1,\ldots, t\}$ such that $T_{1,s}$ and $T_{2,s}$ are a $p$-distinguishing pair;
\item the $6t$ triples consisting of the $2t$ triples $T_{i,s}$ and their cyclic permutations are mutually inequivalent.
\end{enumerate}
Then $G:=H^k$ is a Beauville group for each $k=6t,\ldots, d_2(H)$.
\end{lemma}

\noindent{\sl Proof.} The proof is similar. The first $3t$ coordinates of $a_i, b_i$ and $c_i$ are defined as before, but now those in positions $j=3t+1,\ldots, 6t$ are the coordinates of $a_{3-i}, b_{3-i}$ or $c_{3-i}$ in positions $j-3t$. Each of the six generators $a_i, b_i, c_i$ of $G$ now has at least one $p$-full coordinate for each prime $p$ dividing $|H|$, so any $g\in\Sigma_1^{(p)}\cap\Sigma_2^{(p)}$ satisfies $(1)$ for $i=1$ and $2$, leading to a contradiction as before. \hfill$\square$

\medskip


For the small values of $k$ omitted by these lemmas one can often use the following:

\begin{lemma}\label{strongcoprimeperiods}
Let $H$ be a non-abelian finite simple group with a set of $r\ge 2$ mutually inequivalent generating triples of type $(l,m,n)$, where $l, m$ and $n$ are mutually coprime. Then $G:=H^k$ is a Beauville group for each $k=2,\ldots, 6r$.
\end{lemma}

\noindent{\sl Proof.} If the specified generating triples for $H$ are $(x_j, y_j, z_j)$ for $j=1, \ldots, r$, then one can form $6r$ mutually inequivalent generating triples for $H$ by cyclically permuting each $(x_j, y_j, z_j)$ and each $(z_j^{-1},y_j^{-1},x_j^{-1})$. If $k=2,\ldots, 6r$ one can form two generating triples
\[a_1=(x_1, x_2, \ldots),\quad b_1=(y_1, y_2,\ldots),\quad c_1=(z_1, z_2,\ldots)\]
and
\[a_2=(x_1, y_2, \ldots),\quad b_2=(y_1, z_2,\ldots),\quad c_2=(z_1, x_2,\ldots),\]
for $G$, using $k$ of these triples in the different coordinate positions. Since $l, m$ and $n$ are mutually coprime, any prime $p$ dividing $l$, $m$ or $n$ divides exactly one of them, so if $g\in\Sigma_1^{(p)}$ then $|{\rm supp}(g)\cap\{1, 2\}|=2$ whereas if $g\in\Sigma_2^{(p)}$ then $|{\rm supp}(g)\cap\{1, 2\}|=1$. Thus $\Sigma_1^{(p)}\cap\Sigma_2^{(p)}=\emptyset$, so these two triples form a Beauville structure for $G$. \hfill$\square$

\section{The groups $L_2(q)$}\label{L2q}

We will now apply the results in the preceding section to show that the simple groups $L_2(q)$ satisfy Theorem~\ref{mainthm}.


Let $G=H^k$ where $H=L_2(q)$ with $q=p_0^e\ge 4$ for some prime $p_0$, and with $k=2,\ldots, d_2(H)$. Since $L_2(4)\cong L_2(5)\cong A_5$ and $L_2(9)\cong A_6$, the main theorem of~\cite{Jon13}, which proves Theorem~\ref{mainthm} for $H=A_n$ ($n\ge 5$), allows us to assume that $q=7$ or $8$ or $q\ge 11$. We will use Lemma~\ref{stronglypdist} with $t=1$ to show that $G$ is a Beauville group. However, this lemma does not apply when $k=2$, so we will deal with this case first by a separate argument.

Let $k=2$. Results of Macbeath~\cite{Mac} show that there exist generating triples $(x_i, y_i, z_i)$ $(i=1, 2$) for $H$ of types $(q_1, q_1, p_0)$ and $(q_2, q_2, p_0)$, where $q_1, q_2=(q\pm 1)/d$ and $d=\gcd(2, q-1)$. For instance, in the equation
\[
\Big(\,\begin{matrix}s&1\cr -1&0\end{matrix}\,\Big)
\Big(\,\begin{matrix}0&1\cr -1&t\end{matrix}\,\Big)
=
\Big(\,\begin{matrix}-1&s+t\cr 0&-1\end{matrix}\,\Big)
\]
in $SL_2(q)$ one can choose $s$ and $t$ to be the traces of elements $x_i$ and $y_i$ of order $o(x_i)=o(y_i)=q_i$ in $H$ for $i=1, 2$, with the matrix on the right representing an element $z_i^{-1}=x_iy_i$ of order $p_0$ provided $s+t\ne 0$. Dickson's description of the subgroups of $L_2(q)$ (see~\cite[Ch.~XII]{Dic}) shows that the only maximal subgroups containing elements of orders $q_i$ and $p_0$ are dihedral groups of order $q_i$ when $p_0=2$, impossible here since $x_i$ and $y_i$ do not commute, and stabilisers of points in ${\mathbb P}^1(q)$, also impossible since $z_i$ fixes only $\infty$, which is not fixed by $x_i$. Thus each triple $(x_i, y_i, z_i)$ generates $H$.
Having different types, these two generating triples are mutually inequivalent, so we obtain a generating triple
\[a_1=(x_1, x_2),\quad b_1=(y_1, y_2),\quad c_1=(z_1, z_2)\]
for $G$. Similarly $G$ has a second generating triple
\[a_2=(x_2, z_1),\quad b_2=(y_2, x_1),\quad c_2=(z_2, y_1).\]

The primes $p$ dividing $|H|=q(q^2-1)/d$ are $p_0$ and those dividing $q_1$ or $q_2$. Since $p_0$, $q_1$ and $q_2$ are mutually coprime, if $p=p_0$ then any $g\in\Sigma_i^{(p)}$ has $|{\rm supp}(g)|=2$ or $1$ for $i=1$ or $2$ respectively, so $\Sigma_1^{(p)}\cap\Sigma_2^{(p)}=\emptyset$. If $g\in\Sigma_i^{(p)}$ with $p$ dividing $q_1$, then ${\rm supp}(g)=\{1\}$ or $\{2\}$ for $i=1$ or $2$, while for $p$ dividing $q_2$ it is the other way round. Thus  $\Sigma_1^{(p)}\cap\Sigma_2^{(p)}=\emptyset$ for all $p$, so $H^2$ is a Beauville group. 

\smallskip

Now let $k=3,\ldots, d_2(H)$. Let $(x_i, y_i, z_i)$ ($i=1, 3$) be generating triples for $H$ of types $(q_1, q_1, p_0)$ and $(q_1, q_2, q_2)$ for $H$. For instance, the first could be as above in the case $k=2$;
the second could consist of the images in $H$ of matrices in $SL_2(q)$ of the form
\[
\Big(\,\begin{matrix}u&1\cr -1&0\end{matrix}\,\Big),
\quad
\Big(\,\begin{matrix}v&w\cr 0&v^{-1}\end{matrix}\,\Big)
\quad\hbox{and}\quad
\Big(\,\begin{matrix}-w&-uw-v^{-1}\cr v&uv\end{matrix}\,\Big),
\]
where $u$ is the trace of an element of $H$ of order $q_1$, $v$ generates the multiplicative group of the field ${\mathbb F}_q$, and $w$ is chosen so that $uv-w$ is the trace of an element of $H$ of order $q_2$. Again, it follows from~\cite[Ch.~XII]{Dic} that this triple generates $H$.

We now apply Lemma~\ref{stronglypdist}, with $t=1$, to this pair of triples. Again, we must consider the primes $p$ dividing $|H|$. The generator $z_1$ is $p_0$-full, since the Sylow $p_0$-subgroups of $H$ are elementary abelian, while $x_1, y_1$ and $x_3$ are $p$-full for all primes $p$ dividing $q_1$, and $y_3$ and $z_3$ are $p$-full for primes dividing $q_2$. The hypotheses of Lemma~\ref{stronglypdist} are satisfied, with $t=1$, so $H^k$ is a Beauville group for each $k=3,\ldots, d_2(H)$.

\section{Counting triples}\label{Ct}

For the other families of finite simple groups we shall consider, it is easier to prove the existence of generating triples of various types by means of the following enumerative formula of Frobenius~\cite{Fro} than to exhibit them directly as we did for $H=L_2(q)$.

\begin{prop}\label{Frob}
If $X, Y$ and $Z$ are conjugacy classes in a finite group $H$, then the number $\nu_H(X, Y, Z)$ of triples $(x, y, z)\in X\times Y\times Z$ such that $xyz=1$ is given by
\[\nu_H(X, Y, Z)=\frac{|X|.|Y| .|Z|}{|H|}\sum_{\chi}\frac{\chi(x)\chi(y)\chi(z)}{\chi(1)},\]
where $x\in X$, $y\in Y$, $z\in Z$, and $\chi$ ranges over the irreducible complex characters of $H$.
\end{prop} 

Proofs of this and related results can be found in~\cite[Appendix]{Jon95} or~\cite[Theorem~7.2.1]{Ser}. By summing over all choices of conjugacy classes $X, Y$ and $Z$ of elements of given orders $l, m$ and $n$, one can find the number $\nu_H(l, m, n)$ of triples of type $(l, m, n)$ in $H$.

An equivalent, and sometimes more convenient, version of Proposition~\ref{Frob} is that
\[\nu_H(X, Y, Z)=\frac{|H|^2}{|C_H(x)|.|C_H(y)|.|C_H(z)|}S(X, Y, Z),\]
where $S(X, Y, Z)$ denotes the character sum
\[S(X, Y, Z)=\sum_{\chi}\frac{\chi(x)\chi(y)\chi(z)}{\chi(1)}.\]
In most of the cases we shall consider, one finds that many characters $\chi$ take the value $0$ on at least one of the three classes $X, Y$ and $Z$, so they do not contribute to this sum. Even when a non-principal character $\chi$ does contribute, the values of $|\chi|$ on these classes are often so much smaller than the degree $\chi(1)$ that the sum is dominated by the contribution, equal to $1$, from the principal character, so that $S(X, Y, Z)\approx 1$. These estimates can often show that $S(X, Y, Z)>0$, so that $\nu_H(X, Y, Z)>0$ and hence $H$ has such triples $(x, y, z)$.

We are particularly interested in generating triples, those not contained in any maximal subgroup $M$ of $H$. The number $\phi_H(l, m, n)$ of these of type $(l, m, n)$ clearly satisfies
\[\phi_H(l, m, n)\ge \nu_H(l, m, n)-\sum_M\nu_M(l, m, n),\]
where $M$ ranges over the maximal subgroups of $M$. In particular, if $\nu_M(l,m,n)=0$ for all such $M$ then
\[\phi_H(l, m, n) = \nu_H(l, m, n).\]

\section{The Suzuki groups}\label{Suz}

We will apply this method first to the Suzuki groups $Sz(q)$, using the description of their conjugacy classes, subgroups and characters given by Suzuki in~\cite{Suz} (see also~\cite[\S4.2]{Wil}).

If $q=2^e$ for some odd $e\ge 3$, then the group $H=Sz(q)={}^2B_2(q)$ is a simple group of order $q^2(q^2+1)(q-1)$. It has four conjugacy classes of maximal cyclic subgroups, of mutually coprime orders $4$, $q_1=q+\sqrt{2q}+1$, $q_2=q-\sqrt{2q}+1$ and $q_3=q-1$. (Note that $q_1q_2=q^2+1$.) It follows that any prime $p$ dividing $|H|$ divides precisely one of these four integers, and any non-identity element of $G$ is conjugate to an element of precisely one of these four subgroups. In particular, a pair of generating triples of types $(4, q_2, q_3)$ and $(q_1, q_3, q_3)$, if they exist, will satisfy the conditions of Lemma~\ref{stronglypdist} with $t=1$, and hence prove that $H^k$ is a Beauville group for $k=3,\ldots, d_2(H)$.

To prove that such triples exist, we use Frobenius's formula (Proposition~\ref{Frob}), together with the character table of $H$ in~\cite[\S17]{Suz}. Conjugacy classes $X$, $Y$ and $Z$ of elements of orders $4$, $q_2$ and $q_3$ satisfy $|X|=4(q^2+1)(q-1)$, $|Y|=q^2q_1(q-1)$ and $|Z|=q^2(q^2+1)$, so each choice of such classes gives rise to $4q_1|H|>0$ triples $(x, y, z)\in X\times Y\times Z$ with $xyz=1$. (Only the principal character $\chi$ makes a non-zero contribution to the sum in Proposition~\ref{Frob}, all other characters vanishing on at least one of the chosen classes.) The subgroups of $H$ are described in~\cite[\S 15]{Suz}; no proper subgroup contains elements of orders $q_2$ and $q_3$, so each of these triples generates $H$. A similar argument applies to triples of type $(q_1, q_3, q_3)$: each appropriate choice of conjugacy classes yields $q_2(q+1)|H|$ generating triples, with only the irreducible characters of degrees $1$ and $q^2$ contributing to the sum.

The case $k=2$ is covered by Lemma~\ref{strongcoprimeperiods}: since $|{\rm Aut}\,H|=e|H|$ and $4q_1>4q>e$ there are at least two inequivalent generating triples of type $(4, q_2, q_3)$, and their periods are mutually coprime. Thus the Suzuki groups satisfy Theorem~\ref{mainthm}.


\section{The small Ree groups}\label{Ree}

Similar methods can be applied to the small Ree groups $H=R(q)={}^2G_2(q)$, where $q=3^e$ for some odd $e\ge 3$. (These groups, announced and described by Ree in~\cite{Ree60, Ree61}, are called `small' to distinguish them from the `large' Ree groups ${}^2F_4(2^e)$.) Each of these groups $H$ is simple, of order $q^3(q^3+1)(q-1)$. The Sylow $2$-subgroups of $H$ are elementary abelian, of order $8$; there is a single conjugacy class of involutions, with centralisers isomorphic to $C_2\times L_2(q)$, of order $q(q^2-1)$. The Sylow $3$-subgroups are non-abelian, of order $q^3$ and exponent $9$; there are three conjugacy classes of elements of order $9$, all with centralisers of order $3q$. There are cyclic Hall subgroups $A_i$ of mutually coprime odd orders $q_i:=(q-1)/2$, $(q+1)/4$, $q-\sqrt{3q}+1$ and $q+\sqrt{3q}+1$ for $i=0,\ldots, 3$.
(Note that $q\equiv 3$ mod~$(8)$ and $q^3+1=(q+1)(q^2-q+1)=4q_1q_2q_3$.)

We will again use Lemmas~\ref{stronglypdist} and~\ref{strongcoprimeperiods}, now applied to triples of types $(2, q_1, q_3)$ and $(9, q_0, q_2)$, with six mutually coprime periods. For each prime $p$ dividing $|H|$, one of these triples contains a single $p$-full element, while the elements of the other triple all have orders coprime to $p$. Provided they exist, such triples therefore satisfy the hypotheses of these two lemmas.

One can count triples by applying Frobenius's formula to the character table for $H$ given by Ward in~\cite{War}. In the case of triples of type $(9, q_0, q_2)$, there are respectively  three, $\varphi(q_0)/2$ and $\varphi(q_2)/6$ conjugacy classes of elements of orders $9$, $q_0$ and $q_2$, denoted by $YT^i\;(i=0, \pm 1)$, $R^a$ and $V$ in~\cite{War}, with centralisers of orders $3q$, $q-1$ and $q_2$. Only the principal character makes a non-zero contribution to the character sum in the formula, so the number of triples of this type in $H$ is
\[3.\frac{\varphi(q_0)}{2}.\frac{\varphi(q_2)}{6}.\frac{|H|^2}{3q(q-1)q_2}
=
\frac{1}{12}\varphi(q_0)\varphi(q_2)q^5(q^3+1)(q-1)^2q_3>0.\]

A similar calculation applies to triples of type $(2, q_1, q_3)$. The elements of orders $2$, $q_1$ and $q_3$, denoted by $J$, $S^a$ and $W$ in~\cite{War}, form $1$, $\varphi(q_1)/6$ and $\varphi(q_3)/6$ conjugacy classes, with centralisers of orders $q(q^2-1)$, $q+1$ and $q_3$ respectively. The character sum takes the value
\[1+\frac{q.(-1)^2}{q^3}+2.\frac{q_0.(-1).1}{q_0q_2\sqrt{q/3}}
=
1+\frac{1}{q^2}-\frac{2}{q_2\sqrt{q/3}},
\]
with $\xi_3$ of degree $q^3$, and $\xi_6$ and $\xi_8$ of degree $q_0q_2\sqrt{q/3}$ the only non-principal characters contributing to it. It follows that the number of triples of this type is
\[\frac{\varphi(q_1)}{6}.\frac{\varphi(q_3)}{6}.\frac{|H|^2}{q(q^2-1)(q+1)q_3}
\left( 1+\frac{1}{q^2}-\frac{2}{q_2\sqrt{q/3}}\right)
>0.
\]

The maximal subgroups of $H$ have been determined by Levchuk and Nuzhin~\cite{LN} and Kleidmann~\cite{Kle} (see also~\cite[\S 4.5.3]{Wil}). They are as follows:

\begin{itemize}
\item the point stabilisers in the doubly transitive permutation representation of degree $q^3+1$, or equivalently, the normalisers of the Sylow $3$-subgroups, of order $q^3(q-1)$;
\item the centralisers of involutions, of order $q(q^2-1)$;
\item the normalisers of Hall subgroups $A_i$ ($i=1,2,3$), of order $6q_i$;
\item subgroups isomorphic to $R(3^f)$ where $e/f$ is prime.
\end{itemize}
None of these subgroups contains elements of orders $q_0$ and $q_2$, or of orders $q_1$ and $q_3$, so each triple of type $(9, q_0, q_2)$ or $(2, q_1, q_3)$ generates $H$. Dividing the above numbers of triples by $|{\rm Aut}\,H|=e|H|$ gives the numbers of equivalence classes of triples of these two types. As in the case of the Suzuki groups, Lemmas~\ref{stronglypdist} and~\ref{strongcoprimeperiods} now show that these Ree groups satisfy Theorem~\ref{mainthm}.


\section{The sporadic simple groups}\label{Spo}

In this section we will prove that the $26$ sporadic simple groups satisfy Theorem~\ref{mainthm}. Most of the information we need about these groups can be found in~\cite{ATLAS} or~\cite{Wil}; in particular, we have adopted the notation of~\cite{ATLAS} for conjugacy classes, maximal subgroups, etc. 

\subsection{The Mathieu groups}

It is convenient to deal with the Mathieu groups together, in view of the inclusion relations between them. Each of the Mathieu groups $H=M_n$ ($n=11, 12, 22, 23$ or $24$) has order divisible by at most six primes. This allows us to use Lemma~\ref{stronglypdist} with $t=1$ to deal with the values $k=3,\ldots, d_2(H)$, while Lemma~\ref{strongcoprimeperiods} deals with $k=2$. In each case we can use Frobenius's formula (Proposition~\ref{Frob}), applied to the character tables in~\cite{ATLAS}, to count triples of a given type, and the lists of maximal subgroups in~\cite{ATLAS} to eliminate those which do not generate $H$. Here, as with some of the other sporadic simple groups, we will compute the exact numbers of equivalence classes of generating triples of the required types, even though it is usually sufficient to show that at least one or two exist: this is because it may be useful in other contexts (such as regular maps and hypermaps) to know the number of normal subgroups of a given triangle group with quotient isomorphic to $H$. For instance, by the calculation in \S 9.2.1 there are, up to isomorphism, $18$ orientably regular maps of type $\{5, 19\}$ with orientation-preserving automorphism group isomorphic to Janko's group $J_1$; the Riemann-Hurwitz formula shows that they all have genus $21715$.

\subsubsection {$M_{11}$}

Let $H=M_{11}$, a simple group of order $2^4.3^3.5.11$ and exponent $2^3.3.5.11$ with $|{\rm Out}\,H|=1$. Frobenius's formula and the character table in~\cite{ATLAS} show that there are $2^6.3^3.5^2.11=20|H|=20|{\rm Aut}\,H|$ triples of type $(3, 8, 11)$ in $H$; they all generate $H$ since no maximal subgroup contains elements of orders $8$ and $11$, so they form $20$ equivalence classes.

Similarly, there are $54|{\rm Aut}\,H|$ triples of type $(5, 11, 11)$. The only maximal subgroups containing elements of order $5$ and $11$ are the twelve isomorphic to $L_2(11)$. Applying Frobenius's formula to $L_2(11)$ shows that $8|{\rm Aut}\,H|$ of these triples are contained in such maximal subgroups, so the remainder form $46$ equivalence classes of generating triples.

A pair of generating triples of types $(3, 8, 11)$ and $(5, 11, 11)$ satisfy the hypotheses of Lemma~\ref{stronglypdist}, so $H^k$ is a Beauville group for each $k=3,\ldots, d_2(H)$. Since there are at least two inequivalent generating triples of type $(3, 8, 11)$, Lemma~\ref{strongcoprimeperiods} proves this for $k=2$.

\subsubsection{$M_{12}$}

The group $H=M_{12}$, of order $2^6.3^3.5.11$ with $|{\rm Out}\,H|=2$, has the same exponent as $M_{11}$, so we can use the same method, with triples of the same types $(3, 8, 11)$ and $(5, 11, 11)$, provided they generate $H$. Frobenius's formula shows that $H$ has $96|{\rm Aut}\,H|$ triples of type $(3, 8, 11)$. The only maximal subgroups containing elements of orders $8$ and $11$ are the $24$ isomorphic to $M_{11}$; we have already counted such triples in $M_{11}$, so we find that there are $20|{\rm Aut}\,H|$ non-generating triples, and hence $76$ equivalence classes of generating triples of this type.

Similarly there are $158|{\rm Aut}\,H|$ triples of type $(5, 11, 11)$ in $H$. The only maximal subgroups containing elements of orders $5$ and $11$ are the $24$ isomorphic to $M_{11}$ and the $144$ isomorphic to $L_2(11)$. Having already counted triples of this type in these two groups, we find that there are $54|{\rm Aut}\,H|$ and $4|{\rm Aut}\,H|$ triples contained in maximal subgroups respectively isomorphic to $M_{11}$ and to $L_2(11)$. These two sets of triples are disjoint, since any triple of type $(5, 11, 11)$ contained in $L_2(11)$ generates that group, so there are $58|{\rm Aut}\,H|$ non-generating triples, and hence $100$ equivalence classes of generating triples of this type. From this point, as with the remaining Mathieu groups, the argument follows that for $M_{11}$.

\subsubsection{$M_{22}$}

The group $H=M_{22}$ has order $2^7.3^2.5.7.11$ and exponent $2^3.3.5.7.11$ with $|{\rm Out}\,H|=2$. There are $2^{17}.3^3.5.7.11=1536|{\rm Aut}\,H|$ triples of type $(5, 11, 11)$ in $H$. The only maximal subgroups containing elements of order $11$ are the $2^5.3.7=672$ isomorphic to $L_2(11)$. We have seen that $L_2(11)$ contains $2^5.3.5.11=5280$ such triples, and is generated by each of them, so there are $2^{10}.3^2.5.7.11=4|{\rm Aut}\,H|$ non-generating triples and hence $1532$ equivalence classes of generating triples of this type.

Similarly, there are $216|{\rm Aut}\,H|$ triples of type $(3, 7, 8)$ in $H$, each generating $H$ or contained in a maximal subgroup isomorphic to $AGL_3(2)$. There are $330$ such maximal subgroups, each containing $224$ elements of order $3$ (eight in each of the $28$ two-point stabilisers, isomorphic to $S_4$, in the natural affine action of this group), and $384$ elements of order $7$ ($48$ in each of the eight point stabilisers, isomorphic to $L_3(2)$). Thus there are at most $330\times 224\times 384=32|{\rm Aut}\,H|$ non-generating triples and hence at least $184$ equivalence classes of generating triples of type $(3, 7, 8)$ in $H$.

\subsubsection{$M_{23}$}

The group $H=M_{23}$ has order $2^7.3^2.5.7.11.23$ and exponent $2^3.3.5.7.11.23$ with $|{\rm Out}\,H|=1$. There are $2^{10}.3^2.5.7^2.11^2.23=616|{\rm Aut}\,H|$ triples of type $(3, 8, 23)$. The only maximal subgroups with elements of order $23$ are the normalisers of Sylow $23$-subgroups, of order $11.23$, so each of these triples generates $H$, giving $616$ equivalence classes.

There are $2^{15}.3^3.5.7.11.23^2=17664|{\rm Aut}\,H|$ triples of type $(5, 7, 11)$ in $H$. The only maximal subgroups with elements of orders $5, 7$ and $11$ are the $23$ subgroups isomorphic to $M_{22}$. Now there are $2^{16}.3^4.5.7.11$ triples of type $(5, 7, 11)$ in $M_{22}$, each generating $M_{22}$, so $H$ has $2^{16}.3^4.5.7.11.23 = 4608|{\rm Aut}\,H|$ non-generating triples and hence $13056$ equivalence classes of generating triples of this type.

\subsubsection{$M_{24}$}

The group $H=M_{24}$ has order $2^{10}.3^3.5.7.11.23$ and $|{\rm Out}\,H|=1$, with the same exponent $2^3.3.5.7.11.23$ as $M_{23}$, so we can again use triples of types $(3, 8, 23)$ and $(5, 7, 11)$. There are $2^{11}.3^2.5.7.11.23.769=1538|{\rm Aut}\,H|$ triples of type $(3, 8, 23)$ in $H$ with the generator of order $3$ in the class $3A$, and $2^{14}.3^3.5.7.11^2.17.23=1992|{\rm Aut}\,H|$ with it in class $3B$, giving a total of $3530|{\rm Aut}\,H|$ triples. The maximal subgroups of $H$ have been determined by Choi~\cite{Cho}. The only maximal subgroups with elements of orders $8$ and $23$ are the $24$ isomorphic to $M_{23}$. As we have seen, $M_{23}$ has $2^{10}.3^2.5.7^2.11^2.23$ triples of type $(3, 8, 23)$, each generating $M_{23}$, so there are $2^{13}.3^3.5.7^2.11^2.23=616|{\rm Aut}\,H|$ non-generating triples of this type in $H$, and hence $2914$ equivalence classes of generating triples.

There are $22256|{\rm Aut}\,H|$ triples of type $(5, 7, 11)$ in $H$. The only maximal subgroups with elements of orders $5, 7$ and $11$ are the $24$ isomorphic to $M_{23}$ and the ${24\choose 2}=276$ isomorphic to $M_{22}:2={\rm Aut}\,M_{22}$. These are the stabilisers of points and of unordered pairs of points in the natural representation of $H$. Any triple of type $(5,7,11)$ in ${\rm Aut}\,M_{22}$ must be contained in the subgroup $M_{22}$ of index $2$ fixing two points, so it is contained in two of the point stabilisers. However, no such triple lies in three point stabilisers, since elements of order $11$ fix just two points. We have seen that there are $2^{15}.3^3.5.7.11.23^2=17664|{\rm Aut}\,M_{23}|=17664|{\rm Aut}\,H|/24$ triples of type $(5, 7, 11)$ in $M_{23}$, and $2^{16}.3^4.5.7.11=192|{\rm Aut}\,H|/23$ in $M_{22}$, so the number of generating triples of this type is
\[22256|{\rm Aut}\,H|-17664|{\rm Aut}\,H|+276.\frac{192|{\rm Aut}\,H|}{23}
=6896|{\rm Aut}\,H|,\]
and they form $6896$ equivalence classes.

\subsection {Other small sporadic simple groups}

The method applied to the Mathieu groups can also be applied to some of the other sporadic simple groups $H$, provided they are `small' in the sense that not too many primes divide $|H|$. In these cases, it is sufficient to show that there are inequivalent generating triples satisfying the hypotheses of Lemmas~\ref{stronglypdist} and~\ref{strongcoprimeperiods}. We will deal with these groups in ascending order of magnitude.

\subsubsection{The Janko group $J_1$}

The first Janko group $H=J_1$ has order $2^3.3.5.7.11.19$ and exponent $2.3.5.7.11.19$ with $|{\rm Out}\,H|=1$. We will use triples of types $(2,5,19)$ and $(3,7,11)$. There are $2^4.3^2.5.7.11.19=18|{\rm Aut}\,H|$ triples of type $(2, 5, 19)$ in $H$. They are all generating triples, since no maximal subgroup (classified by Janko~\cite{Jan66}) contains elements of orders $5$ and $19$, so they form $18$ equivalence classes.
Similarly there are $2^5.3.5.7.11.19^2=76|{\rm Aut}\,H|$ triples of type $(3, 7, 11)$ in $H$. No maximal subgroup contains elements of orders $7$ and $11$, so these form $76$ equivalence classes of generating triples.

\subsubsection{The Janko group $J_2$}

The second Janko group $H=J_2$ has order $2^7.3^3.5^2.7$ and exponent $2^3.3.5.7$ with $|{\rm Out}\,H|=2$. We will use triples of types $(3, 5, 7)$ and $(7, 7, 8)$. There are $2^{16}.3^4.5^2.7=768|{\rm Aut}\,H|$ triples of type $(7, 7, 8)$ in $H$, with only the principal character and that of degree $225$ contributing to the character sum. The only triples of this type contained in proper subgroups are those which generate one of the $100$ maximal subgroups isomorphic to $U_3(3)$. This group contains $2^{12}.3^3.7$ such triples, so there are $2^{14}.3^3.5^2.7=64|{\rm Aut}\,H|$  non-generating triples in $H$, and hence there are $704$ equivalence classes of generating triples of this type.

Any triple of type $(3, 5, 7)$ generates $H$ since no maximal subgroup has elements of orders $5$ and $7$. For simplicity we will consider only those triples for which the element of order $3$ is in the ${\rm Aut}\,H$-invariant conjugacy class $3B$, consisting of those with a centraliser of order $36$, and the element of order $5$ is in class $5C$ or $5D$, transposed by ${\rm Aut}\,H$, with a centraliser of order $50$. We find that there are $49$ equivalence classes of such triples, with only the principal character and that of degree $288$ appearing in the character sum.

\subsubsection{The Higman-Sims group $HS$}

The Higman-Sims group $H=HS$ has order $2^9.3^2.5^3.7.11$ and exponent $2^3.3.5.7.11$ with $|{\rm Out}\,H|=2$. We will use triples of types $(7, 11, 11)$ and $(8, 11, 15)$. The only maximal subgroups containing elements of order $11$ are isomorphic to $M_{11}$ or $M_{22}$; neither of these groups has elements of order $15$, and only $M_{22}$ has elements of order $7$. There are $104760|{\rm Aut}\,H|$ triples of type $(7, 11, 11)$ in $H$, with only the principal character and that of degree $3200$ contributing to the character sum. There are $100$ subgroups isomorphic to $M_{22}$, in total containing (and generated by) $2048|{\rm Aut}\,H|$ such triples, so there are $102712$ equivalence classes of triples of this type. There are $50496|{\rm Aut}\,H|$ triples of type $(8, 11, 15)$, with only the principal character and that of degree $175$ contributing to the character sum. Each of these triples generates $H$, so they form $50496$ equivalence classes.




\subsubsection{The McLaughlin group $M^cL$}

The McLaughlin group $H=M^cL$ has order $2^7.3^6.5^3.7.11$ and exponent $2^3.3^2.5.7.11$ with $|{\rm Out}\,H|=2$. We will use triples of types $(5, 7, 11)$ and $(8, 9, 9)$. In counting triples of type $(5, 7, 11)$ we will consider only those for which the element of  order $5$ is in the conjugacy class $5A$, with centralisers of order $750$. The only maximal subgroups with elements of orders $7$ and $11$ are those isomorphic to $M_{22}$; the permutation character of $H$ on the cosets of such a subgroup, as given in~\cite{ATLAS}, takes the value $0$ on elements of $5A$,  so these elements do not lie in such subgroups. Thus each of these triples generates $H$. There are $2^{12}.3^{11}.5^3.7.11=3888|{\rm Aut}\,H|$ of them, with only the principal character contributing to the character sum, so they form $3888$ equivalence classes.

There are $2^{13}.3^6.5^6.7^2.11^2=2^5.5^3.7.11|{\rm Aut}\,H|=308000|{\rm Aut}\,H|$ triples of type $(8, 9, 9)$ in $H$, again with only the principal character contributing to the character sum. The only maximal subgroups containing elements of orders $8$ and $9$ are those isomorphic to $U_4(3)$, $3_+^{1+4}$:$2S_5$ or $3^4$:$M_{10}$; in the last two cases, any triple of this type must lie in the subgroup $3_+^{1+4}$:$2A_5$ or $3^4$:$A_6$ of index $2$, whereas this has no elements of order $8$. Thus any non-generating triple of this type must lie in a maximal subgroup isomorphic to $U_4(3)$. There are $275$ of these, any two of them intersecting in a subgroup isomorphic to $3^4$:$A_8$ or $L_3(4)$; these have no elements of orders $8$ or $9$ respectively, so each such triple lies in a unique maximal subgroup. Since $U_4(3)$ has $2^{14}.3^7.5.7.23$ triples of type $(8, 9, 9)$ (with only the principal character and the two of degree $35$ contributing to the character sum), there are $2^{14}.3^7.5^3.7.11.23=2^7.5.7|{\rm Aut}\,H|=4480|{\rm Aut}\,H|$ non-generating triples in $H$, and hence $303584|{\rm Aut}\,H|$ generating triples, forming $303584$ equivalence classes.

\medskip

One can apply this method to a few other sporadic simple groups, such as $J_3$, $He$, $Ru$, $O$'$N$ and $Co_3$, but it is more efficient to deal with them by using a separate method.

\subsection{The larger sporadic simple groups}

In the case of the larger sporadic simple groups $H$, and in particular those of order divisible by more than six primes, Lemma~\ref{stronglypdist} may not be applicable, since it may be impossible to find a pair of generating triples satisfying its hypotheses for each of the primes dividing $|H|$. Instead, one can use Lemma~\ref{pdist} with $t>1$ to deal with values $k\ge 6t$, and Lemma~\ref{strongcoprimeperiods} to deal with small values of $k$. In each case, we find two mutually coprime integers $q_1$ and $q_2$ dividing $|H|$ (usually the largest two prime divisors) such that no maximal subgroup of $H$ contains elements of orders $q_1$ and $q_2$. (The maximal subgroups have been completely classified in all cases except that of the Monster simple group $M$, and in this case sufficient is known about the maximal subgroups to justify this claim.) It follows that for any $l$, a triple of type $(l, q_1, q_2)$ in $H$ must generate the whole group. Provided they exist, we can therefore use pairs of triples of types $(l, q_1, q_2)$ for various $l$ to distinguish the primes dividing $|H|$. By Lemma~\ref{pdist} this will deal with large values of $k$, and one of these types, with $l$ coprime to $q_1$ and $q_2$, will also deal with small $k$ by Lemma~\ref{strongcoprimeperiods}, provided there are sufficiently many inequivalent triples of that type.

As before, we can use Frobenius's formula (Proposition~\ref{Frob}) to show that suitable triples exist. In all the relevant cases, inspection of the character table of $H$ shows that the character sum $S(X, Y, Z)$ is very close to $1$, with the modulus $|\chi(g)|$ of each non-principal character $\chi$ of $H$ for $g$ in each of the chosen conjugacy classes $X, Y$ and $Z$ significantly smaller than the degree $\chi(1)$. Thus $\nu_H(X, Y, Z)>0$, so suitable generating triples exist, and moreover in sufficient numbers for Lemma~\ref{strongcoprimeperiods} to deal with the small values of $k$ not covered by Lemma~\ref{pdist}.

We will give a detailed explanation of this method  as it applies to the largest and the smallest of the remaining sporadic simple groups, namely the Monster group and the third Janko group. For the other groups, in order to omit tedious and repetitive numerical details we will simply refer to Table~1, which lists the groups in ascending order of magnitude, together with their exponents and the choices of $q_1$ and $q_2$.

\begin{table}[htb]
\centering
\begin{tabular}{|c|c|c|c|}
\hline
$H$&name&exponent&$q_1, q_2$\\
\hline\hline
$J_3$&Janko&$2^3.3^2.5.17.19$&$17, 19$\\
\hline
$He$&Held&$2^5.3.5.7.17$&$7, 17$\\
\hline
$Ru$&Rudvalis&$2^3.3.5.7.13.29$&$13, 29$\\
\hline
$Suz$&Suzuki&$2^3.3^2.5.7.11.13$&$11, 13$\\
\hline
$O$'$N$&O'Nan&$2^4.3.5.7.11.19.31$&$17, 31$\\
\hline
$Co_3$&Conway&$2^4.3^2.5.7.11.23$&$9, 23$\\
\hline
$Co_2$&Conway&$2^4.3^2.5.7.11.23$&$9, 23$\\
\hline
$Fi_{22}$&Fischer&$2^4.3^2.5.7.11.13$&$11, 13$\\
\hline
$HN$&Harada-Norton&$2^3.3^2.5^2.7.11.19$&$11, 19$\\
\hline
$Ly$&Lyons&$2^3.3^2.5^2.7.11.31.37.67$&$37, 67$\\
\hline
$Th$&Thompson&$2^3.3^3.5.7.13.19.31$&$19, 31$\\
\hline
$Fi_{23}$&Fischer&$2^4.3^3.5.7.11.13.17.23$&$17, 23$\\
\hline
$Co_1$&Conway&$2^4.3^2.5.7.11.13.23$&$13, 23$\\
\hline
$J_4$&Janko&$2^4.3.5.7.11.23.29.37.43$&$37, 43$\\
\hline
$Fi'_{24}$&Fischer&$2^4.3^3.5.7.11.13.17.23.29$&$23, 29$\\
\hline
$B$&Baby monster&$2^5.3^3.5^2.7.11.13.17.19.23.31.47$&$31, 47$\\
\hline
$M$&Monster&$2^5.3^3.5^2.7.11.13.17.19.23.29.31.41.47.59.71$&$59, 71$\\
 \hline
\end{tabular}
\caption{The larger sporadic simple groups}
\end{table}

\subsubsection{The Monster group $M$}

The Monster group $H=M$ has order $2^{46}.3^{20}.5^9.7^6.11^2.13^3.17.19.23.29.31.41.47.59.71$ and exponent $2^5.3^3.5^2.7.11.13.17.19.23.29.31.41.47.59.71$, with $|{\rm Out}\,H|=1$. We will apply the method described above, with $q_1=59$ and $q_2=71$. The maximal subgroups of $M$ are not yet completely known, but it is known that the only maximal subgroups of order divisible by $59$ or $71$ are those isomorphic to $L_2(59)$ or $L_2(71)$ (see~\cite[Table~5.6]{Wil}); thus no proper subgroup has order divisible by both $59$ and $71$, so a triple of type $(l, 59, 71)$ for any $l$ must generate $H$. Provided they exist, we can therefore use a pair of triples of types $(59, 59, 71)$ and $(59, 71, 71)$ to distinguish the primes $59$ and $71$, and six pairs of triples of types $(l, 59, 71)$, where $l=11, 13, 19, 23, 25, 27, 29, 31, 32, 41, 47, 119\;(=7.17)$, for the remaining thirteen primes $p<59$ dividing $|H|$. By Lemma~\ref{pdist} this will deal with values $k\ge 42$, and one of the latter types, with mutually coprime periods, will cover smaller $k$ by Lemma~\ref{strongcoprimeperiods}. Inspection of the character table of $H$ in~\cite{ATLAS} shows that for each of the above types one can choose appropriate conjugacy classes $X, Y$ and $Z$ so that the character sum $S(X, Y, Z)$ is very close to $1$. This shows that such generating triples exist, and in sufficient numbers for Lemma~\ref{strongcoprimeperiods} to deal with the cases $k=2,\ldots, 41$.

In fact, all but $14$ of the $193$ non-principal irreducible characters $\chi$ vanish on elements of order $59$ or $71$, so they do not contribute to any of the character sums, while the non-vanishing characters, of  degrees $\chi(1)\ge 8980616927734375$, all take values $\chi(g)=\pm 1$ or $(-1\pm i\sqrt{59})/2$ on those classes. The character values on the other classes needed are also much smaller than these degrees: for instance, if $g$ has order $l=11$ then $|\chi(g)|\le 190$ for all $\chi$. Since there are two conjugacy classes each of elements of orders $59$ and $71$, and one of order $11$, with centralisers of orders $59, 71$ and $1045440$ respectively, it follows that the number of equivalence classes of generating triples of type $(11, 59, 71)$ is approximately
\[4\times\frac{|H|}{59.71.1045440}\approx 7\cdot 38\times 10^{44}.\]
The estimates for the other types required are even larger, since there are more elements of order $l$ in the cases where $l>11$.

\subsubsection{The Janko group $J_3$}

The third Janko group $H=J_3$ has order $2^6.3^5.5.17.19$ and exponent $2^3.3^2.5.17.19$ with $|{\rm Out}\,H|=2$. Although there are only five primes dividing $|H|$, it is more convenient to apply the method used for $M$ than that used for the Mathieu groups. We will take $q_1=17$ and $q_2=19$. The only maximal subgroups of $H$ with elements of order $19$ are those isomorphic to $L_2(19)$; this group has no elements of order $17$, so for any $l$, a triple of type $(l, 17, 19)$ must generate $H$. Provided they exist, one can therefore use a pair of triples of types $(17, 17, 19)$ and $(17, 19, 19)$ to distinguish the primes $17$ and $19$, a pair of types $(8, 17, 19)$ and $(9, 17, 19)$ for the primes $p=2$ and $3$, and a pair of types $(5, 17, 19)$ and (for instance) $(4,17,19)$ for $p=5$. As in the case of $M$, the character values guarantee the existence of all the required triples: the only non-principal irreducible characters $\chi$ not vanishing on elements of orders $17$ or $19$, and hence contributing to any of the relevant character sums $S(X, Y, Z)$,  are one of degree $324$, two of degree $1215$ and three of degree $1920$; these satisfy $|\chi(g)|\le 4$ for all elements $g$ of order at least $4$, so each $S(X, Y, Z)$ is very close to $1$. By Lemma~\ref{pdist} this deals with the values $k=18, \ldots, d_2(H)$. Frobenius's formula shows that there are, for instance, $21312$ equivalence classes of triples of type $(5, 17, 19)$, more than enough for Lemma~\ref{strongcoprimeperiods} to deal with the values $k=2,\ldots, 17$.

\medskip

Table~1 shows how the same method can be applied to the remaining $15$ sporadic simple groups. Thus each of the $26$ sporadic simple groups satisfies Theorem~\ref{mainthm}, and the proof is complete.


\section{The groups $L_3(q)$ and $U_3(q)$}\label{L3qU3q}

The method used for the larger sporadic simple groups can also be applied to other families of simple groups, such as the groups $L_3(q)$ and $U_3(q)$.

\subsection{$L_3(q)$}
The $3$-dimensional projective special linear group $H=L_3(q)$ is simple for every prime power $q$. Since $L_3(2)\cong L_2(7)$ we may assume that $q>2$. The group $H$ contains elements $h_1$ and $h_2$ of orders
\[q_1=\frac{q^2+q+1}{d}\quad {\rm and}\quad q_2=\frac{q^2-1}{d}\]
where $d=\gcd(q^2+q+1, 3)=\gcd(q-1,3)$. The maximal subgroups of $H$ were determined for even $q$ by Hartley~\cite{Har}, and for odd $q$ first by Mitchell~\cite{Mit} and then by Bloom~\cite{Blo} (see also~\cite{KlLi}): there are none containing elements of orders $q_1$ and $q_2$ (see also Berecky's classification~\cite{Ber} of subgroups of classical groups containing Singer cycles). It follows that for any $l$, a triple of type $(l, q_1, q_2)$ must generate $H$. The generic character table for $H$ given by Simpson and Frame in~\cite{SF} shows that $\chi(h_1)\chi(h_2)=0$ for every non-principal irreducible character $\chi$ of $H$ except the Steinberg character, of degree $q^3$; for this character we have $\chi(h_1)=1=-\chi(h_2)$, and $|\chi(h)|\le q$ for all $h\ne 1$. Thus $S(X,Y,Z)>0$ where $X$ is any non-identity conjugacy class, and $Y$ and $Z$ are classes of elements of orders $q_1$ and $q_2$, so generating triples of type $(l, q_1, q_2)$ exist whenever $H$ contains elements of order $l\ne 1$.

Now
\[|H|=\frac{1}{d}q^3(q^2+q+1)(q+1)(q-1)^2,\]
so each prime divisor $p$ of $|H|$ divides at least one of $q$, $q^2+q+1$, $q+1$ and $q-1$. These four integers are mutually coprime, except that $\gcd(q^2+q+1, q-1)=3$ if $q\equiv 1$ mod~$(3)$, and $\gcd(q+1, q-1)=2$ if $q$ is odd. It follows from this that $q_1$ and $q_2$ are always mutually coprime. A pair of generating triples of types $(q_1, q_1, q_2)$ and $(q_1, q_2, q_2)$ can therefore be used to strongly distinguish all the primes dividing $q_1$ or $q_2$. This deals with every prime $p$ dividing $|H|$, except the prime $p=p_0$ dividing $q$, and the prime $p=3$ if $q\equiv 4$ or $7$ mod~$(9)$ (so that $q_1$ and $q_2$ are both coprime to $3$, whereas $|H|$ is not). The exponent $q_0$ of a  Sylow $p_0$-subgroup of $H$ is $p_0$ if $p_0>2$, and $4$ if $p_0=2$, so a pair of triples of types $(q_0, q_1, q_2)$ and $(3, q_1, q_2)$ will strongly distinguish $p_0$ if $p_0\ne 3$, and will strongly distinguish $3$ if $q\equiv 4$ or $7$ mod~$(9)$; if $p_0=3$, so that $q\not\equiv 4$ or $7$ mod~$(9)$, we can instead replace the second triple with one of type $(2 , q_1, q_2)$ to strongly distinguish $p_0$.  In either case these two pairs of triples satisfy the hypotheses of Lemma~\ref{stronglypdist} with $t=2$, so they deal with the values $k=6,\ldots, d_2(H)$. 

There are at least two inequivalent generating triples of type $(p_0, q_1, q_2)$: for instance, one can choose the first element to be central or non-central in the unique Sylow $p_0$-subgroup containing it. These triples have mutually coprime periods, so Lemma~\ref{strongcoprimeperiods}(1) deals with $k=2,\ldots, 12$, and hence the groups $H=L_3(q)$ satisfy Theorem~\ref{mainthm2}.

\subsection{$U_3(q)$}

The proof for the unitary groups $H=U_3(q)$ is essentially the same. These groups have order
\[|H|=\frac{1}{d}q^3(q^3+1)(q-1)=\frac{1}{d}q^3(q^2-q+1)(q+1)(q-1)\]
where $d=\gcd(q+1,3)$, and they  are simple for all prime powers $q>2$. They contain elements $h_1$ and $h_2$ of mutually coprime orders
\[q_1=\frac{q^2-q+1}{d}\quad {\rm and}\quad q_2=\frac{q^2-1}{d},\]
all satisfying $H=\langle h_1, h_2\rangle$ (see~\cite{KlLi} or~\cite[\S 3.10.9]{Wil} for the maximal subgroups of $H$, or~\cite{Ber} for subgroups containing a Singer cycle $h_1$). As with $L_3(q)$, the generic character table for $H$ in~\cite{SF} shows that generating triples of type $(l, q_1, q_2)$ exist whenever $H$ contains elements of order $l\ne 1$. A pair of generating triples of types $(q_1, q_1, q_2)$ and $(q_1, q_2, q_2)$ strongly distinguish each prime dividing $|H|$, apart from the prime $p_0$ dividing $q$, and the prime $3$ if $q\equiv 2$ or $5$ mod~$(9)$. The exponent $q_0$ of a  Sylow $p_0$-subgroup is as before, so for these primes one can use a pair of triples of types $(q_0, q_1, q_2)$ and either $(3, q_1, q_2)$ or $(2 , q_1, q_2)$ as $q\equiv 2, 5$ mod~$(9)$ or not. Lemma~\ref{stronglypdist}, with $t=2$, then deals with the values $k=6,\ldots, d_2(H)$, and again Lemma~\ref{strongcoprimeperiods}(1), applied to triples of type $(p_0, q_1, q_2)$, deals with small $k$. This completes the proof of Theorem~\ref{mainthm2}.


\end{document}